\documentclass[11pt]{article}

\usepackage{amsmath, amsfonts, amsthm, amssymb}
\usepackage{graphicx}

 \DeclareMathOperator{\SL}{SL}
 
\begin{document}
\def \docdate{August 3, 2007}
\pagestyle{myheadings}
\markboth{qq}{Continued Fractions \hfill \docdate \hfill M. Pavlovskaia\hfill }

\title{Continued Fraction Expansions of Matrix Eigenvectors}
\author{Maria Pavlovskaia \\
maria\_pavlovskaia@hmc.edu \\
\date \docdate}
\date{\today}

\maketitle

\begin{abstract}
We examine various properties of the continued fraction expansions of
matrix eigenvector slopes of matrices
 from the $\SL$(2, $\mathbb Z$) group. We calculate the average period length,
maximum period length, average period sum, maximum period sum
and the distributions of 1s 2s and 3s in the periods versus the radius
of the Ball within which the matrices are located. We also prove that
the periods of continued fraction expansions from the real irrational roots of  
$x^{2}+px+q=0$ 
 are always palindromes.

\end{abstract}

\section{Introduction}

Any real number $x$ can be expressed as a continued fraction of this form:
\begin{equation*}
x = a_{0}+\cfrac{1}{a_{1}+  
\cfrac{1}{a_{2}+
\cfrac{\dotsb}{a_{n-1}+
\cfrac{1}{a_{n}+\dotsb} }}}
\end{equation*}
where $a_{0} \in \mathbb Z$ and $a_{i} \in \mathbb Z^{+}$ for all $i
\geq 1$. For simplicity, we write  $ x = [a_{0},a_{1},\dotsb, a_{n}, \dotsb]$.
We will be looking at quadratic irrationalities, or real roots
of the equation $ax^{2}+bx+c=0$
with integer $a, b, c$
. When a root is irrational, 
its continued fraction expansion is infinite, and 
always periodic past a certain index.
Any rational has a finite continued fraction expansion.

In the Continued Fractions brochure by V.I. Arnold
\cite{Arnold2}
he proposes this problem: 
"Let us examine matrices 
$\begin{bmatrix}
a&b\\ c&d
\end{bmatrix}
$, with integral $a,b,c,d$ and the determinant equal to 1.
1. We choose from these the ones that define a hyperbolic rotation.
There is a finite number of matrices whose coefficients are not too large,
ie. $a^{2}+b^{2}+c^{2}+d^{2} \leq N^{2}$.
For each such matrix there exists an eigenvector $y=\lambda x$,
for which $\lambda$ is a quadratic irrationality and therefor the 
continued fraction expansion of $\lambda$ is periodic. 
Take this period and calculate how many ones, twos, threes, etc. there are in
its period, and then average this for all matrices $\begin{bmatrix}
a&b\\ c&d
\end{bmatrix}$. That is to say, take the number of ones in each period
and divide by the number of elements in the period for each matrix.
Hypothesis: this relation will approach a Gaussian distribution as
$N$ approaches infinity." (paraphrased into English from Russian.
See also \cite{Arnold3}). M.Avdeeva and B.Bykovski solved
another problem posed in this broshure regarding the Gaussian
distribution of period elements of continued fractions \cite{Bykovski}.

We have written a program that generates matrices $M \in \SL(2, \mathbb Z)$ , such
that M = $\begin{bmatrix}
k&l\\ m&n
\end{bmatrix}
$ is within a Ball or radius $r$ around the origin  for all 
$r \leq 1000$. That is, it generates matrices with the coefficients
$k, l, m, n$ such that 
$k^{2}+l^{2}+m^{2}+n^{2} \leq r^{2}$
and $kn-lm=1$
.  We then calculate the continued fraction expansions of
the matrix eigenvector slopes
 and find some statistics about their periods. 
 Because any quadratic irrationality can be expressed
 as the slope of an eigenvector of such a matrix,
 all possible quadratic irrationalities are included in 
 our statistics as $r$ approaches infinity. 

I'd like to extend a special thank you to Oleg Karpenkov for valuable help with the revision
process.

\section{Methods}
Our program generates all matrices $M \in \SL(2, \mathbb Z)$
in a ball of radius $r$ around the origin, such
that M = $\begin{bmatrix}
k&l\\ m&n
\end{bmatrix}
$ and $k \geq 0$
 for all 
$r\leq 1000$ using the algorithm desribed in 2.1. We find one eigenvector 
$\begin{bmatrix}
1\\ y_{1}
\end{bmatrix}$
for each matrix $M$ using $y_{1}= \frac{n-k+\sqrt{(k+n)^{2}-4(kn-lm)}}{2l}$.
If $y_{1}$ is real, we find the period of its continued fraction expansion
using a variation of the Euclidean Algorithm described in 2.2.
If $y_{1}$ is rational, we say that its period is 0.
We do not make any calculations 
for the second eigenvector
$\begin{bmatrix}
1\\ y_{2}
\end{bmatrix}$ of M because by Lemma 1, 
$y_{2}$ has the same
coefficients in its period as $y_{1}$, but in reverse order. 
We count each matrix with coefficient $k \geq 1$ twice, because it has the same
eigenvectors as the matrix $M^{'}=$  
$\begin{bmatrix}
-k&-l\\ -m&-n
\end{bmatrix}$ which we had not accounted for.

\subsection{Matrix Generating Algorithm}
We generate all matrices $M \in \SL(2, \mathbb Z)$
in a ball of radius $r$ around the origin
such
that M = $\begin{bmatrix}
k&l\\ m&n
\end{bmatrix}
$ and $k \geq 0$. 
We do this by generating all $k$ between $0$ and $r$.
Then for each $k$ we generate all $m$ between 
$-\sqrt{r^{2}-k^{2}}$ and $\sqrt{r^{2}-k^{2}}$.
For each $k,m$ pair, we find a solution to the Diaphontane equation
$kn-lm=1$ for $l,n$. This way we obtain one matrix
M = $\begin{bmatrix}
k&l\\ m&n
\end{bmatrix}
$ for each $k,m$ pair. To obtain the rest of the matrices,
we find $M_{a}=\begin{bmatrix}
k&l+ak\\ m&n+am 
\end{bmatrix}
$ for $a=0,1,2,...$ and $a=-1,-2,-3...$
such that $k^{2}+(l+ak)^{2}+m^{2}+(n+am)^{2} \leq r^{2}$.

\subsection{Computing the Continued Fractions of Quadratic Irrationalities}
We wish to compute the continued fraction expansion of 
$x_{0}=\frac{u_{0}+v_{0}\sqrt{w_{0}}}
{z_{0}}$ for $u_{0}$, $v_{0}$, $w_{0}$, $z_{0} \in \mathbb Z$ 
using the Euclidean Algorithm. 
Let $a_{0}$ be the integer part of $x_{0}$ and let $x_{0}^{'}$ 
be the fractional part of $x_{0}$. Then let $x_{1}=\frac{1}{x_{o^{'}}}$
and $a_{1}$ be the integer part of $x_{1}$. We continue this process
finding $x_{0}, x_{1}, \dotsb, x_{n}$ until $x_{n} = x_{i}$ for some
$i < n$. Then the continued fraction expansion of 
$x$ is $[a_{0}, a_{1}, \dotsb,
a_{n}]$ and the period of the continued fraction expansion is $[a_{i}, 
a_{i+1}, \dotsb, a_{n-1}]$.

Notice that $x_{1}=\frac{u_{1}+v_{1}\sqrt{w_{0}}}
{z_{1}}$ for $u_{1}, v_{1}, z_{1} \in \mathbb Z$
that are easily expressed in terms of $u_{0}, v_{0}, z_{0}$ and $a_{0}$.
V. I. Arnold conjectured that if we begin with $x_{0}$ that is a real root of
$x^{2}+px+q$ for $p, q \in \mathbb Z$, then the $u_{i}$ and $z_{i}$ obtained
through this algorithm are always divisible by $v_{i}$. In other words, we are
always able to cancel out $v_{i}$ in this special case. This conjecture held
in all of our calculations. Furthermore, we found that for the general case
where $x_{0}$ is a real root of $ax^{2}+bx+c$ for $a, b, c \in \mathbb Z$
we are also able to cancel out large numbers from the numerator and
denominator by finding their greatest common divisors. We noticed that 
after canceling, the coefficients
$u_{i}, v_{i}, w_{i}$, and $z_{i}$ grow very little throughout the calculation
of any continued fraction expansion. We compared our algorithm with the same
algorithm that does not cancel out the common divisors and found
that our algorithm has a significantly lower big O running time.
Without canceling the coefficients grow exponentially
throughout the calculation, while with
canceling they are bounded.

We also discovered a better algorithm that we have not yet 
implemented.
The idea is to find at every $i$-th step the equation 
$A_{i}x^{2}+B_{i}x+C_{i}$ that has the root $x_{i}$
described above. Then 
if we let $a_{i}$ be the integer part of $x_{i}$,
the equation for $x_{i+1}$ has coefficients
$A_{i+1} = A_{i}a_{i}^{2}+B_{i}a_{i}+C_{i}$
, $B_{i+1} = B_{i}+2A_{i}a_{i}$ 
and $C_{i+1} = A_{i}$
. Hinchin 
\cite{Hinchin} proved that the coefficients $A_{i}, B_{i}$
and $C_{i}$ are bounded with the bounds $C_{i}=A_{i-1}<
2|A_{0}x_{0}|+|A_{0}|+|B_{0}|$ and $B_{i}$ can be
defined in terms of $A_{i}$ and $C_{i}$ by  $B_{i}^{2}
= B_{0}^{2}-4A_{0}C_{0}+4A_{i}C_{i}$.

\section{Results}
\subsection{Average Period Length vs. Radius}
We computed the average period length of the 
continued fraction expansions of the
slopes of
matrix eigenvectors for matrices 
within a ball of radius
$r$ around the origin for $r\leq1000$. We did this by
summing all the period lengths for matrices within a given
radius and dividing by the total number of matrices with
real eigenvectors and real eigenvector slopes within that radius.
We plotted the average period length versus the radius as
shown in Fig. 1. We found that the average period length
grows as $\lg{r}$ and can be very precisely approximated by
$-6.11+1.92\lg{r}$ for $r\leq1000$.

\begin{figure}[htbp]
 \centering
\includegraphics{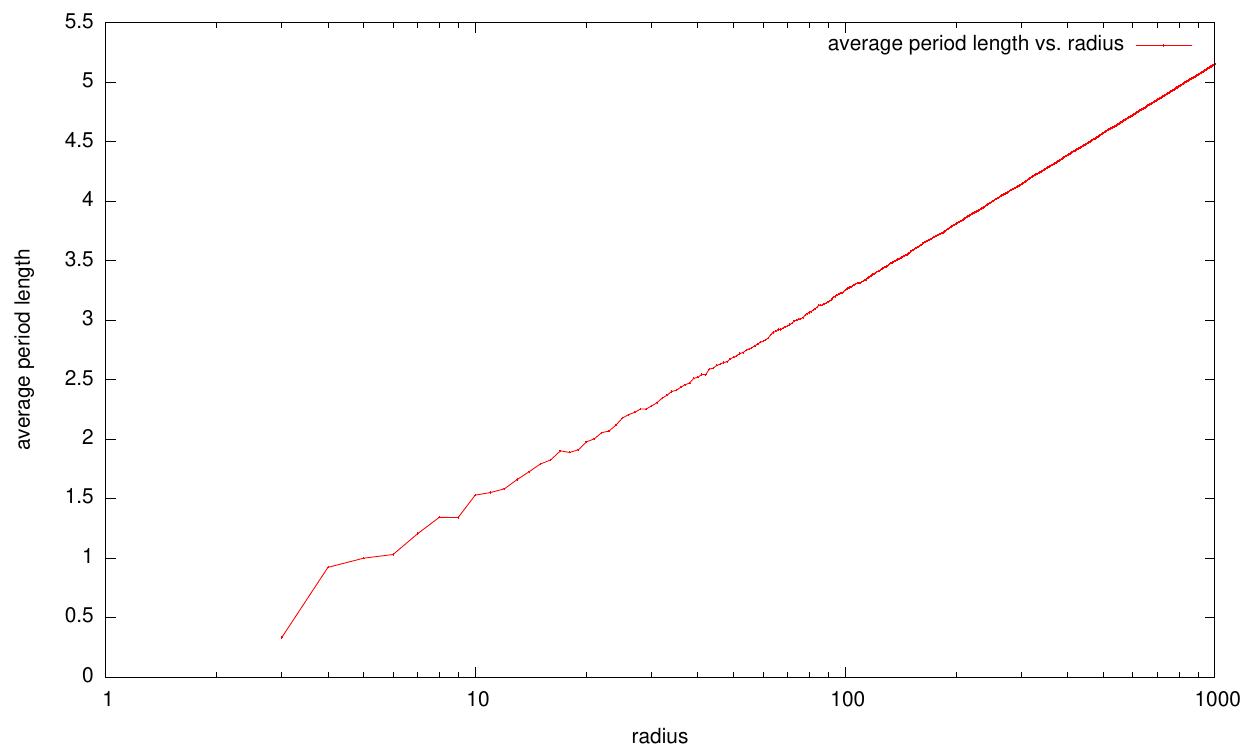}
 \caption{Average period length of 
 continued fractions associated with matrix eigenvectors
 for matrices within a ball of radius r for $r\leq1000$.
It can be very precisely approximated by
$-6.11+1.92\lg{r}$. }
 \label{fig:avelen}
\end{figure}

\subsection{Maximum Period Length vs. Radius}
We computed the maximum period length of the 
continued fraction expansions of the
slopes of
matrix eigenvectors for matrices 
within a ball of radius
$r$ around the origin for $r\leq1000$. 
We plotted the maximum period length versus the radius as
shown in Fig. 2. We found that the maximum period length
appears to grow as $\lg{r}$. 
For example, for a radius of $r=100$, the maximum period length
is 8. It occurs for the matrix 
 $M = \begin{bmatrix}
8&21\\ -29&-76
 \end{bmatrix}
 $, the continued fraction expansion of whose eigenvector
 has the period [1, 1, 1, 1, 1, 1, 1, 2].


\begin{figure}[htbp]
 \centering
\includegraphics{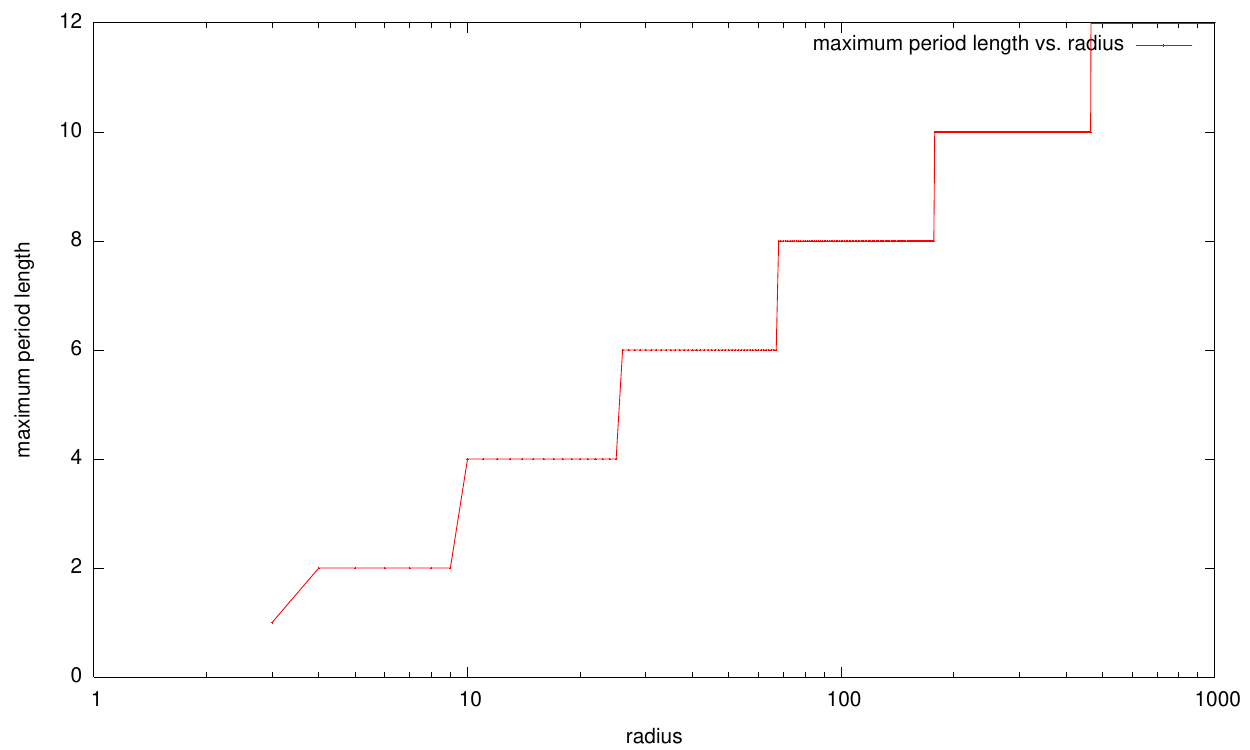}
 \caption{Maximum period length of 
 continued fractions associated with matrix eigenvectors
 for matrices within a ball of radius r. The maximum period length
grows as $\lg{r}$.}
 \label{fig:maxlen}
\end{figure}

\subsection{Average Sum of Period Elements vs. Radius}
We computed the average sum of period elements of the 
continued fraction expansions of the slopes of
matrix eigenvectors for matrices 
within a ball of radius
$r$ around the origin for $r\leq1000$. We did this by
summing all the elements in the periods for matrices within a given
radius and dividing by the total number of matrices with
real eigenvectors and real eigenvector slopes within that radius.
We plotted the average period sum versus the radius as
shown in Fig. 3. We found that the average period sum 
grows faster than $\lg{r}$.

\begin{figure}[htbp]
 \centering
\includegraphics{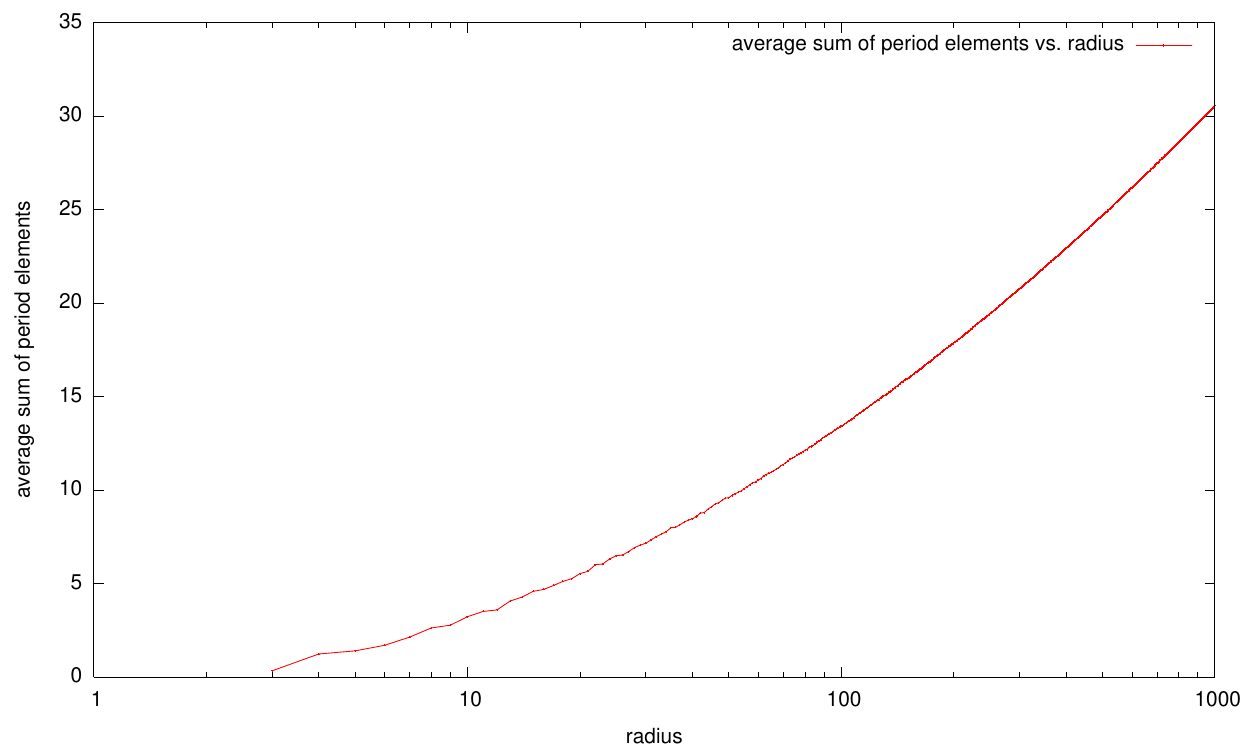}
 \caption{Average sum of period elements of 
 continued fractions associated with matrix eigenvectors
 for matrices within a ball of radius r. The average period sum 
grows faster than $\lg{r}$.}
 \label{fig:avesum}
\end{figure}

\subsection{Average of the Quantity (Period Sum / Period Length) vs. Radius}
We computed the average 
quantity (sum of period elements / period length)
of the periods of the continued fraction expansions of the slopes of
matrix eigenvectors for matrices 
within a ball of radius
$r$ around the origin for $r\leq1000$. We did this by
summing the quantities $Q_{i}$ 
for matrices within a given
radius and dividing by the total number of matrices with
real eigenvectors and real eigenvector slopes within that radius.
Here $Q_{i}$ is the sum of the period elements divided by the length
of the period for the $i$-th continued fraction. 
We plotted the average of the $Q_{i}$ versus the radius as
shown in Fig. 4. We found that the average of the quantity
(period sum/period length)
grows as $\lg{r}$ and can be very precisely approximated by 
$-2.70+3.71\lg{r}$ for $r\leq1000$.

\begin{figure}[htbp]
 \centering
\includegraphics{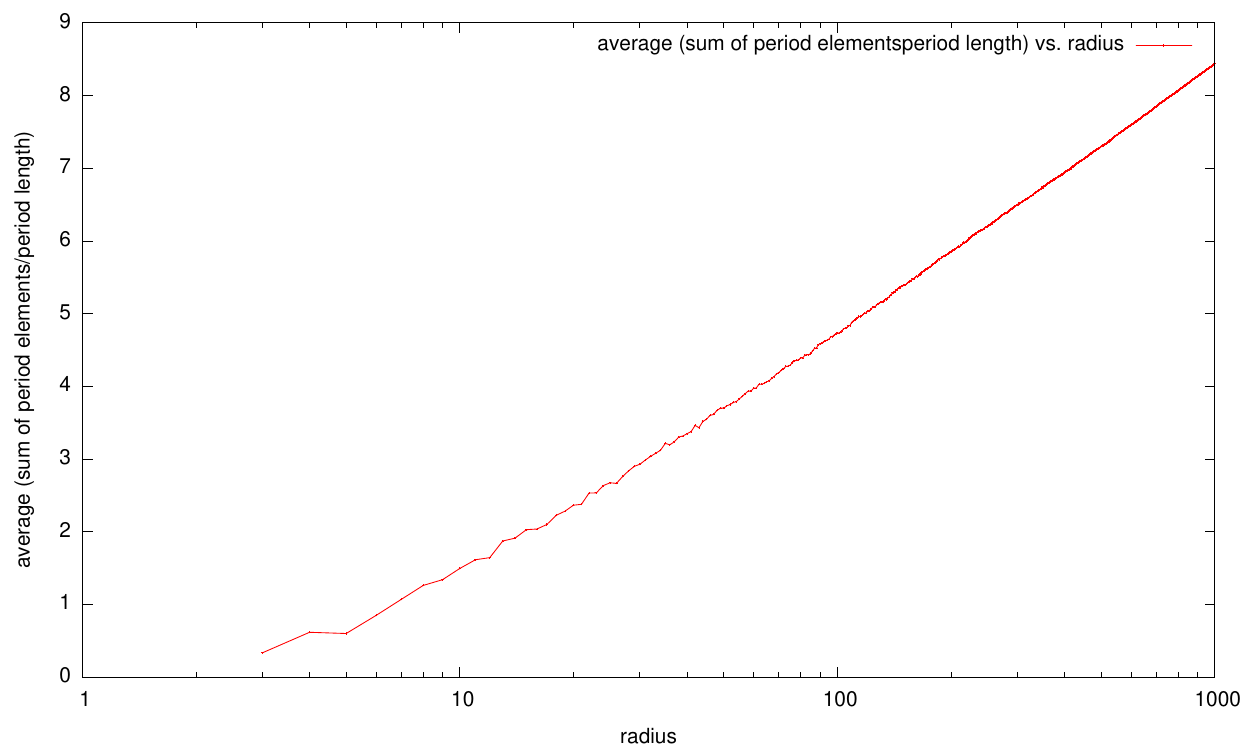}
 \caption{Average of (period sum/period length) for 
 continued fractions associated with matrix eigenvectors
 for matrices within a ball of radius r. It can be
 very precisely approximated by 
$-2.70+3.71\lg{r}$ for $r\leq1000$}
 \label{fig:aveavesum}
\end{figure}

\subsection{Maximum Sum of Period Elements vs. Radius}
We computed the maximum sum of the period elements 
of the continued fraction expansions of the slopes of
matrix eigenvectors for matrices 
within a ball of radius
$r$ around the origin for $r\leq1000$. 
We found that for a radius $r\geq5$ the maximum sum of the period
elements is $r-2$ and occurs for the matrix
 $ M = \begin{bmatrix}
 0&1\\ -1&1-r
 \end{bmatrix}
 $. The eigenvector $e_{1}=$ 
 $\begin{bmatrix}
 1\\ \lambda_{1}
 \end{bmatrix}
 $
 of $M$
 has the slope $\lambda_{1}=$ 
 $\frac{1-r+\sqrt{(1-r)^{2}-4}}{2}$
 , which has the
 continued fraction expansion $[-1, 1, r-3, \dotsb]$ 
 with the period $[1, r-3]$.

\subsection{Appearance of 1s 2s and 3s vs. Radius}
The Gauss-Kuzmin distribution gives the probability distribution
of the occurrence of a given integer in the periods of the continued
fraction expansions of arbitrary real numbers. 
\cite{WikiGauss}. 
Bykovski and Avdeeva proved that this is also 
true for arbitrary quadratic irrationalities
\cite{Bykovski}.
The percentage of 1s, 2s and 3s in 
the periods of the fractions that we calculated should follow
the Gauss-Kuzhmin distribution
\begin{equation*}
\Pr(K=k)=-\log_2\left(1-\frac{1}{(k+1)^{2}}\right).
\end{equation*}
So we would expect that the percentage of 1s would tend 
to 0.415, the percentage of 2s would tend to 0.169 and the 
percentage of 3s would tend to 0.093.
It is evident from Fig. 5, 6 and 7 that the radius is too small
for us to see this distribution.  

\begin{figure}[htbp]
 \centering
\includegraphics{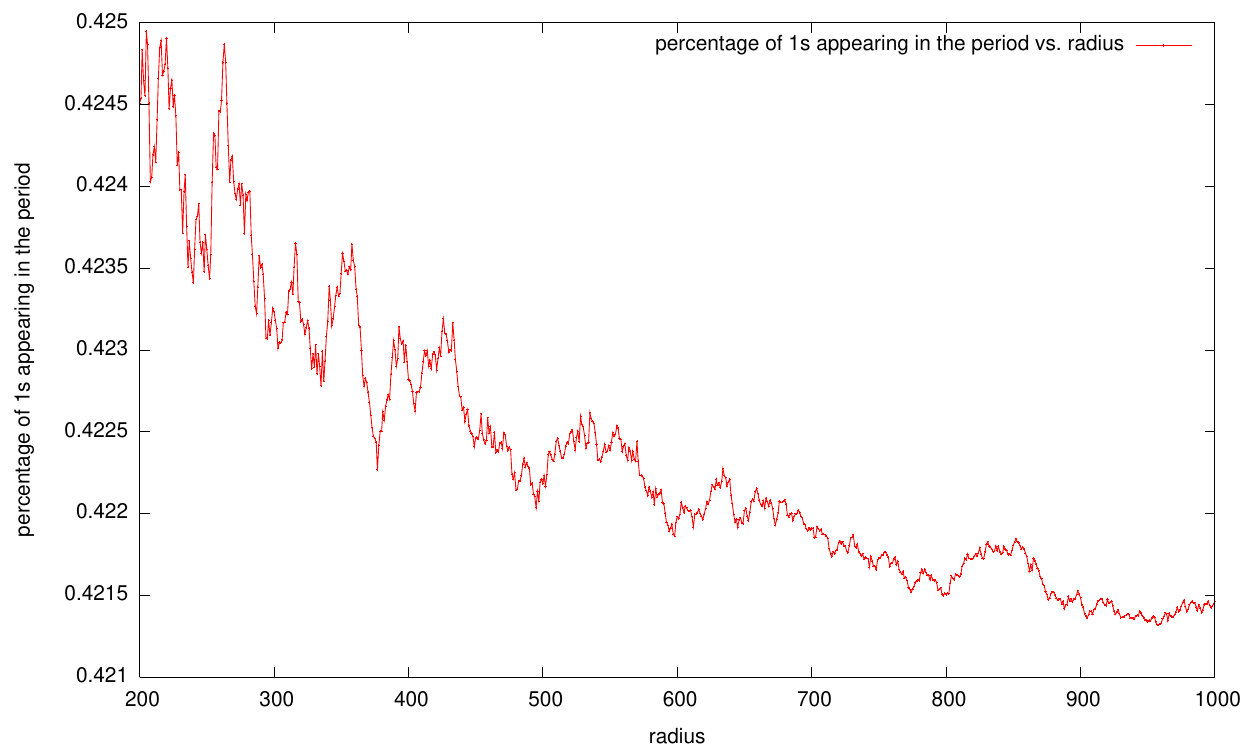}
 \caption{Percentage of 1s in the periods of 
 continued fractions associated with matrix eigenvectors
 for matrices within a ball of radius r
 for $r\leq1000$. The Gauss-Kuzhmin distribution
 requires that the percentage should tend to 0.415 
 as $r$ tends to infinity.}
 \label{fig:aveavesum}
\end{figure}

\begin{figure}[htbp]
 \centering
\includegraphics{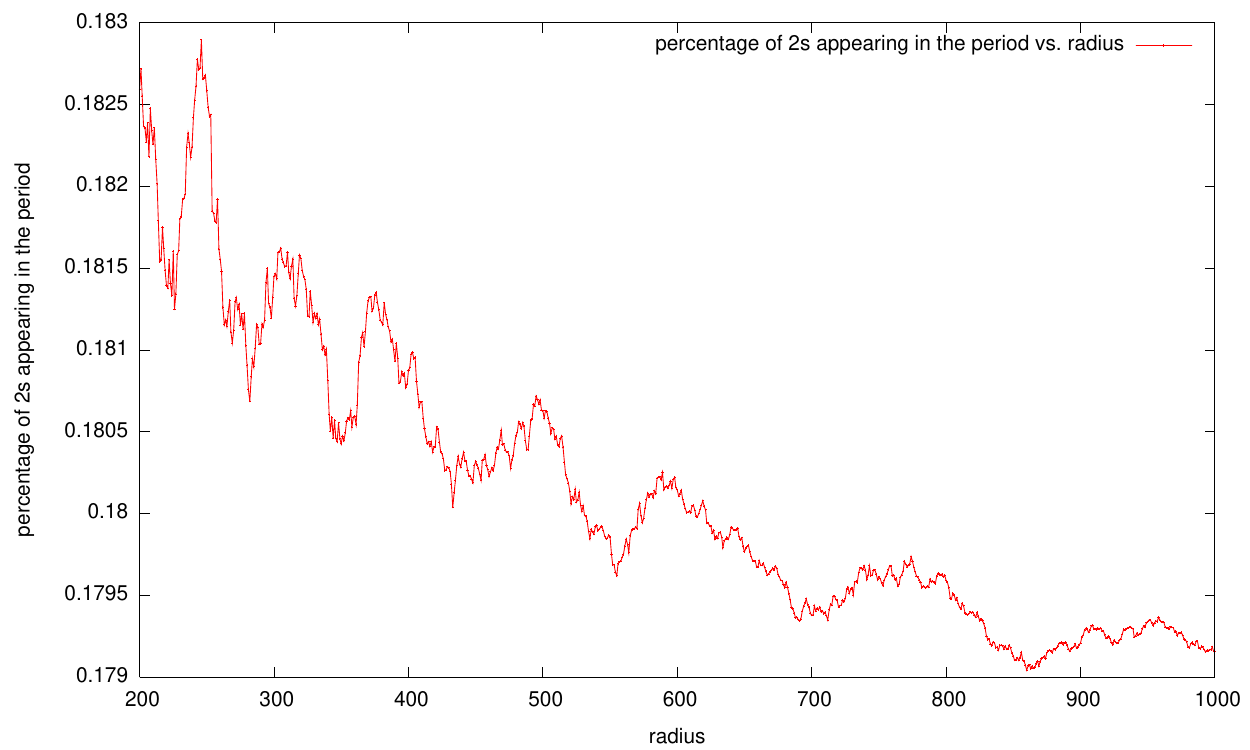}
 \caption{Percentage of 2s in the periods of 
 continued fractions associated with matrix eigenvectors
 for matrices within a ball of radius r
 for $r\leq1000$.The Gauss-Kuzhmin distribution
 requires that the percentage should tend to 0.169 
 as $r$ tends to infinity.}
 \label{fig:aveavesum}
\end{figure}

\begin{figure}[htbp]
 \centering
\includegraphics{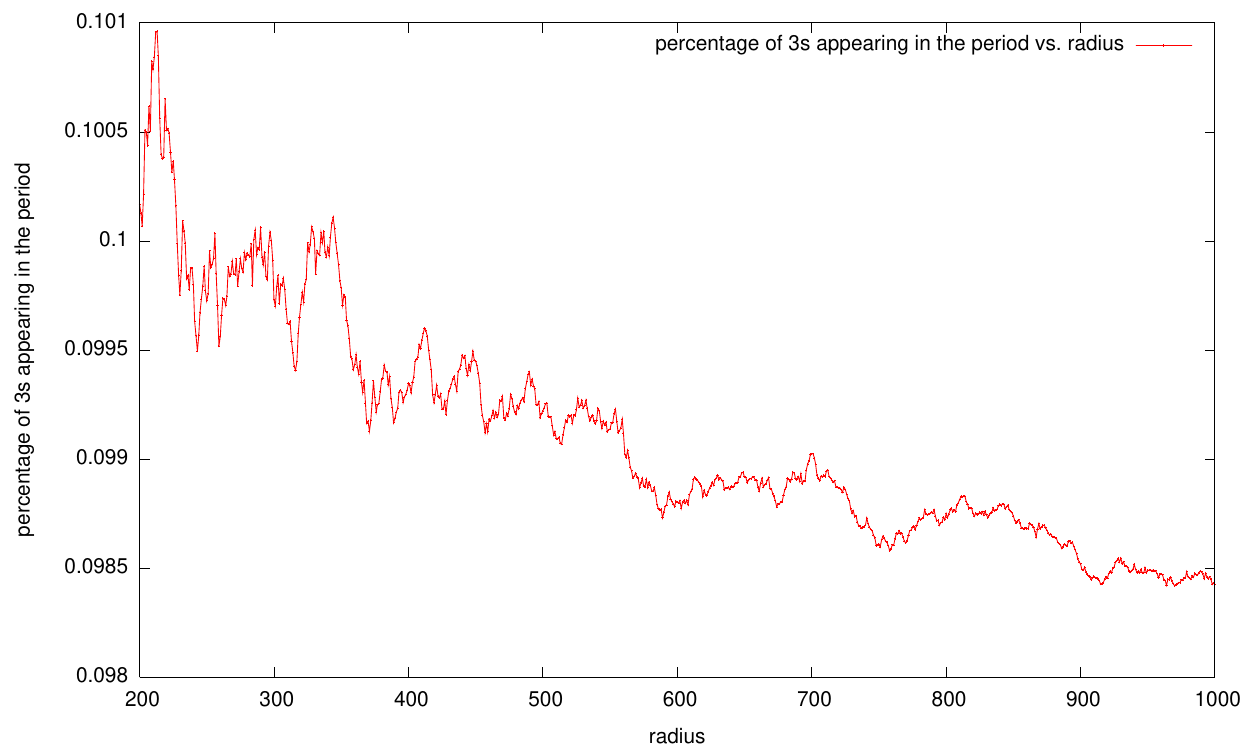}
 \caption{Percentage of 3s in the periods of 
 continued fractions associated with matrix eigenvectors
 for matrices within a ball of radius r
 for $r\leq1000$.The Gauss-Kuzhmin distribution
 requires that the percentage should tend to 0.093 
 as $r$ tends to infinity.}
 \label{fig:aveavesum}
\end{figure}

\section{Palindromes Proof}
\newtheorem{theo}{Theorem}
\newtheorem{lem}{Lemma}  
\newtheorem{prop}{Proposition}
\newtheorem{rem}{Remark}
\begin{theo} Let $x_{1}$ be a real root of $x^{2}+px+q$ for $p, q \in \mathbb Z$.
And let $[a_{0}, a_{1}, \dotsb, a_{n-1}, a_{n}]$ be the period of the continued
fraction of $x_{1}$. Then the period is a palindrome in the sense that 
if we look at the period as a cycle, the cycle is symmetric about some point
(which might be an element of the cycle or be between two elements of the
cycle).
This conjecture was stated by V. I. Arnold 
\cite{Arnold}
 and proven by 
Francesca Aicardi~\cite{Fran}.
Here is an alternate proof that we found independently. 

\begin{proof}
This proof rests on two big ideas. First, if $x_{1}$ and $x_{2}$ are 
real roots of
\begin{equation}
x^{2}+px+q
\end{equation}
for $p, q \in \mathbb Z$, they both have the same period. This was proven by V. I. Arnold~\cite{Arnold}.
Second, Proposition 1 states that if $y_{1}$ and $y_{2}$ are real roots
of any quadratic equation with integer coefficitents
, they have the same coefficients in their period, but in reverse
order. Because $x_{1}$ and $x_{2}$ are roots of (1),
they must have the same period and the period of $x_{1}$ must be the reverse
of the period of $x_{2}$. Consequently, the period of $x_{1}$ must be a palindrome.
\end{proof}
\end{theo}

\begin{prop} Let $y_{1}$ and $y_{2}$ be real roots of
$ay^{2}+by+c$
for $a, b, c \in \mathbb Z$
. If we look at the period
as a cycle, then the period of $y_{1}$ has the same coefficients
as the period of $y_{2}$, but in reverse order.

\begin{proof}
Lemma 1 states that if we let $M \in
\SL(2, \mathbb Z)$ be a matrix with eigenvectors 
$\begin{bmatrix}
1\\ e_{1}
\end{bmatrix}$
and 
$\begin{bmatrix}
1\\ e_{2}
\end{bmatrix}$,
then 
the period of $e_{1}$, $[a_{1}, a_{2}, ..., a_{n}]$ is the reverse 
of the period of $e_{2}$, $[a_{n}, ... a_{2}, a_{1}]$. 
Lemma 2 states that we can construct a matrix $M \in
\SL(2, \mathbb Z)$ such that $\begin{bmatrix}
1\\ y_{1}
\end{bmatrix}$
and  $\begin{bmatrix}
1\\ y_{2}
\end{bmatrix}$
are the eigenvectors of $M$.
It follows that the period of $y_{1}$ is the reverse of the period of $y_{2}$. 

\end{proof}

\end{prop}
\begin{lem}
Let $M$ be a hyperbolic matrix in $
\SL(2, \mathbb Z)$ with eigenvectors 
$\begin{bmatrix}
1\\ e_{1}
\end{bmatrix}$
and 
$\begin{bmatrix}
1\\ e_{2}
\end{bmatrix}$.
Then 
the period of $e_{1}$, $[a_{1}, a_{2}, ..., a_{n}]$ is the reverse 
of the period of $e_{2}$, $[a_{n}, ... a_{2}, a_{1}]$.
\begin{proof}
The proof is based on the geometric interpretation of a continued fraction
as the boundary of the convex hulls of integral points in the angles formed by two lines.
The convex hull formed between the two eigenvectors of
a matrix is invariant under a linear transformation of that matrix. The boundary
is shifted by the matrix along itself. The region of the period on the boundary
is shifted to another region of the period on the boundary.   

Let $\lambda_{1}$ and
$\lambda_{2}$ be eigenvalues of M.
Because $M$ is hyperbolic and
has a determinant of 1, we can assume WLOG that $|\lambda_{1}| > 1$ and $|\lambda_{2}|<1$.
. Therefore, by applying $M$ finitely many times, we can shift the period
of $e_{1}$ onto the period of $e_{2}$. Consequently, their periods must
be the reverse of one-another.

\end{proof}

\end{lem}
\begin{lem}
Let $y_{1}$ and $y_{2}$ be real irrational roots of
$ay^{2}+by+c$ for $b$, $c \in \mathbb Z$. Then
we can construct a matrix $M \in
\SL(2, \mathbb Z)$ such that $\begin{bmatrix}
1\\ y_{1}
\end{bmatrix}$
and  $\begin{bmatrix}
1\\ y_{2}
\end{bmatrix}$
are the eigenvectors of $M$.
\begin{proof}
The continued fraction expansion of $y_{1}$ has some period. It is known that there exists
a matrix $M \in \SL(2, \mathbb Z)$ which shifts one period of $y_{1}$ into the next
one 
and  $\begin{bmatrix}
1\\ y_{1}
\end{bmatrix}$
is an eigenvector of $M$. We wish to shows that 
$\begin{bmatrix}
1\\ y_{2}
\end{bmatrix}$
is the other eigenvector of $M$.
Let $\lambda_{1}$ and $\lambda_{2}$ be eigenvalues of $M$.
Let  $u=\alpha+\beta\sqrt{\gamma}$ 
where $\gamma$ is the discriminant of the characteristic
equation of $M$ and $\alpha$ and $\beta,$ are rational.
We define conjugation for $u$
 as $\overline{u}$
$=\alpha-\beta\sqrt{\gamma}$. 
Then $y_{1}$ satisfies $[M-\lambda_{1}I]$
$\begin{bmatrix}
1\\ y_{1}
\end{bmatrix}$=0.
Consequently it is also true that 
$[\overline{M}-\overline{\lambda_{1}I}]$
$\begin{bmatrix}
\overline{1}\\ \overline{y_{1}}
\end{bmatrix}=0$.
Thus $[M-\lambda_{2}I]$
$\begin{bmatrix}
1\\ y_{2}
\end{bmatrix}=0$, and thus 
$\begin{bmatrix}
1\\ y_{2}
\end{bmatrix}$
is an eigenvector of $M$, as we wanted.

\end{proof}

\end{lem}

\begin{rem}
In all of our calculations the greater root of $x^{2}+px+q$
the period was always either exactly a palindrome or
would be a palindrome with its first or last element removed.
However, this is not always true for the periods of the lesser root. 
We also noticed that for $|p|, |q| \leq 100$ the period of the first
root is never shifted by more than 3 places from that of the second root.

\end{rem}

\bibliographystyle{abbrv}
\bibliography{paper}  

\end{document}